\documentclass[11pt]{article}
\usepackage{amsmath}
\usepackage{mathrsfs}
\usepackage{amsfonts}

\usepackage{amsfonts, amsmath, amssymb, amssymb}
\usepackage{amssymb,amsfonts,amsmath,color,
latexsym, epsfig,cite, psfrag,eepic,colordvi}
\usepackage{amscd,graphics}
\usepackage{lineno}

\newcommand{\qed}{\hfill $\Box $}
\newcommand{\pf}{\noindent {\bf Proof.} }

\newtheorem{theorem}{Theorem}[section]
\newtheorem{lemma}[theorem]{Lemma}
\newtheorem{coro}[theorem]{Corollary}

\newtheorem{conjecture}[theorem]{Conjecture}

\newtheorem{obse}[theorem]{Observation}

\begin{document}

\title{A characterization on orientations of graphs avoiding given lists on out-degrees\thanks{Supported by the National Natural
Science Foundation of China under grant No.12271425}}

\author{
 Xinxin Ma and Hongliang Lu \footnote{Corresponding author: luhongliang215@sina.com}\\School of Mathematics and Statistics\\
Xi'an Jiaotong University\\
Xi'an, Shaanxi 710049, China\\
\smallskip\\
}
\date{}


\date{}

\maketitle

\begin{abstract}
Let $G$ be a graph and $F:V(G)\to2^N$ be a set function. The graph $G$ is said to be \emph{F-avoiding} if there exists an orientation $O$ of $G$ such that $d^+_O(v)\notin F(v)$ for every $v\in V(G)$, where $d^+_O(v)$ denotes the out-degree of $v$ in the directed graph $G$ with respect to $O$.  In this paper, we give a Tutte-type good characterization to decide the $F$-avoiding problem when  for every $v\in V(G)$, $|F(v)|\leq \frac{1}{2}(d_G(v)+1)$ and $F(v)$ contains no two consecutive integers. Our proof also gives a simple polynomial algorithm to find a desired orientation. As a corollary, we prove the following result: if for every $v\in V(G)$,  $|F(v)|\leq \frac{1}{2}(d_G(v)+1)$ and $F(v)$ contains no two consecutive integers, then $G$ is $F$-avoiding. This partly answers a problem proposed by Akbari et. al.(2020)

\end{abstract}

\section{Introduction}
Let $G$ be a  graph without loops and with vertex set $V(G)$ and edge set $E(G)$. Let $e(G):=|E(G)|$. For $v\in V(G)$, we use $d_G(v)$ to denote the degree of a vertex $v$ in $G$. For $S\subseteq V(G)$, let $G[S]$ denote the vertex induced subgraph induced by $S$. Given two  integers $s,t$ such that $s\leq t$, $\{s,s+1,\ldots,t\}$ is denoted by $[s,t]$. An orientation of $G$ is an assignment of a direction to each edge of $G$. For a vertex $v\in V(G)$, we denote by $d^+_O(v)$ the out-degree of $v$ under the orientation $O$ of the edges of $G$.

 Because of its fruitful applications, an orientation with specified properties has been extensively studied\cite{NM,JD}. Frank and Gy\'arf\'as \cite{AA} proved that for a graph $G$ and two mappings $a,b: V(G) \rightarrow \mathbb{N}$ with $a(v)\leqslant b(v)$ for every vertex $v$, $G$ has an orientation such that $$a(v)\leqslant d^+(v)\leqslant b(v)$$ for every vertex $v$, if and only if for any $U\subseteq V(G)$, $$\sum_{v\in U}a(v)-d(U)\leqslant |E(G[U])|\leqslant \sum_{v\in U}b(v).$$ where $d(U)$ is the number of edges connecting $U$ and $V(G)\setminus U$. Borowiecki, Grytczuk and Pil\'{s}niak\cite{MJ} discovered a beautiful fact that every graph
admits an orientation of its edges such that the outdegrees of any two adjacent vertices
are different. Such orientations can be interpreted as graph colorings and are now known
as proper orientations. 

Let $F:V(G)\to2^N$. A graph is said to be \emph{F-avoiding} if there exists an orientation $O$ of $G$ such that $d^+_O(v)\notin F(v)$ for every $v\in V(G)$. In this case, we say that $O$ avoids $F$ or $O$ is an  \emph{F-avoiding} orientation. Conversely, for $H:V(G)\to2^N$, we call that  an orientation $O$ of $G$ an  \emph{H-orientation}, if $d^+_O(v)\in H(v)$ for every $v\in V(G)$. 
 In this paper,  we may always assume that for every $v\in V(G)$, $\max H(v)\leq d_G(v)$ and  $\min H(v)\geq 0$.

 Akbari et al. \cite{SM} proved that for a graph $G$ and $F:V(G)\to2^N$, if
$$|F(v)|\leq \frac{d_G(v)}{4}$$
for every vertex $v$, then $G$ is \emph{F-avoiding}. Furthermore, they proposed that the following conjecture.
\begin{conjecture}[Akbari et al,\cite{SM}]\label{AK}
	For a graph $G$, and a mapping  $F:V(G)\to2^N$, if
	$$|F(v)|\leq \frac{1}{2}(d_G(v)-1)$$
	for every $v\in V(G)$, then $G$ is  \emph{F-avoiding}.
\end{conjecture}
In the same paper, Akbari et al. \cite{SM} proved that Conjecture \ref{AK} holds for bipartite graphs. Bradshaw et.al. made some research progress on Conjecture \ref{AK}: they \cite{BCMMW} proved if $|F(v)|\leq |d_G(v)|/3$ for every $v\in V(G)$, then $G$ is $F$-avoiding.

In order to settle this conjecture, in this paper, we can make some progress with some constraints on the mapping $F$. Our proof is also to give polynomial algorithm to find such orientation.
In this proof, we use the technique of  \emph{parity trace} to describe the existence of the suitable orientation $O$ of $G$ in the constraints of the mapping $H:V(G)\to2^N$, where $H$ has the property:
\begin{align}\label{main-eq}
i\notin H(v)\ \mbox{ implies}\  i+1\in H(v),\ \mbox{for every}\  i\in \mathbb{N},\ v\in V(G).
\end{align}
\medskip
Let us call the pair $(G,H)$ \emph{dense} if $H$ satisfies (\ref{main-eq}), and $G$ is connected. \emph{Parity trace} is a specific permutation of all the cut vertices of $G$, the details of which are described in Section 3.

To formulate main result, we need  to introduce the following notions.
Let $(G,H)$ be dense, $V^1_0$ denote the set of odd vertices, and $V^2_0$ the set of even vertices. Suppose $(x_1,...,x_k)$ is an order of $V(G)\setminus (V^1_0\cup V^2_0)$. $(x_1,...,x_k)$ will be called a \emph{parity trace} if it satisfies the following statement.


\begin{itemize}
\item [(i)] Let $V^1_0:=\{x\ |\ H(x) \mbox{ is odd}\}$ and $V^2_0:=\{x\ |\ H(x) \mbox{ is even}\}$. For $i=1,...k$ define recursively the numbers $l_i, u_i$, and the sets $V^1_i$ and $V^2_i$, for which (\ref{parity}) above holds;
\item [(ii)] Suppose $V^1_{i-1}$ and $V^2_{i-1}$ have already been defined. Let $l_i$ be the number of those components $C$ of $G-x_i$ for which $C\subseteq V^1_{i-1} \cup V^2_{i-1}$ and satisfy $|C\cap  V^1_{i-1}|+ e(C)\not\equiv d_G(x_i,C) \pmod 2$.

 Let $u_i:=d_G(x)-t_i$, where $t_i$ is the number of those components $C$ of $G-x_i$ for which $C\subseteq V^1_{i-1} \cup V^2_{i-1}$ and $|C\cap  V^1_{i-1}|\neq e(C)\pmod 2$. Moreover, the following statement holds.
 \begin{align}\label{parity}
\mbox{All elements in }	[l_i,u_i]\cap H(x_i) \mbox{ have  the  same parity.}
\end{align}

\item [(iii)]  If $[l_i,u_i]\cap H(x_i)$ is odd, then let $V^1_{i}:=V^1_{i-1}\cup x_i$, and $V^2_{i}:=V^2_{i-1}$; if it is even, then $V^1_{i}:=V^1_{i-1}$, and $V^2_{i}:=V^2_{i-1}\cup x_i$. There is no other cases according to (\ref{parity}). ($[l_i,u_i]\cap H(x_i)=\emptyset$ is possible, then we are free to consider it to be either odd or even.)
\end{itemize}
If $(x_1,...x_k)$ is a parity trace, then let $V^1:=V^1_k, V^2:=V^2_k$. Clearly, $V(G)=V^1\cup V^2, V^1\cap V^2=\emptyset$.

In this paper, we give a Tutte-type characterization for $H$-orientation problem when $(G,H)$ is dense.

\begin{theorem}\label{main}
Let $G$ be a graph and let $H:V(G)\rightarrow 2^N$ such that $(G,H)$ is \emph{dense}. There exists an $H$-orientation of $G$ if and only  if there is no parity trace such that $|V^1|\neq e(G)\pmod 2$.
\end{theorem}
 We mainly refer to the idea for finding parity factors in \cite{Guan}, which also appears in \cite{AS}. By Theorem \ref{main}, we may give a partial solution to Conjecture \ref{AK}. In Section 2, we characterize the existence of \emph{H-orientation} for \emph{dense pairs} with 2-connected $G$; this characterization is used in Section 3 to obtain the characterization for connected graphs.
\section{2-Connected Graphs}
First we prove a simple statement concerning arbitrary \emph{dense pairs}. Before proceeding the main arguments, we put several terminologies that are used in next proof. For a given \emph{dense pair} $(G,H)$ and an orientation $O$ of $G$, we call $v\in V(G)$ \emph{feasible} if $d^+_O(v)\in H(v)$.
Let $S\subseteq V(G)$, and let $e(S)$ denote the number of edges whose two vertices both belong to $C$, i.e., $e(S):=|E(G[S])|$. We denote by $d_G(a_i,S)$ the number of edges joining $a_i$ to some vertex of $S$, and $d^+_G(a_i,S)$ the out-degree of $a_i$ to $S$. For $a\in V(G)$, we denote $\delta(a)$ the set of edges incident with $a$. Let $O$ be an  orientation  of $G$ and $P(x,y)$ be a path (or cycle when $x=y$) on $G$ connecting vertices $x$ and $y$. For $ e\in E(G)$, denote $O_e$ the orientation on $e$ based on $O$ and $-O_e$ the reverse orientation with respect to $O_e$.  We define the \emph{symmetry difference} between $O$ and $P(x,y)$ by $O\bigtriangleup P(x,y)$, that is,

$$\hat{O}=O\bigtriangleup P(x,y):=\begin{cases}
	O_e, & \text{if $e\notin P(x,y)$;}\\
	-O_e, &\text{if $e\in P(x,y)$.}
\end{cases}$$

\begin{lemma}\label{prop1}
Let $u\in V(G)$.	If $(G,H)$ is dense,  then there exists an orientation $O$ of $G$ such that $v$ is feasible for all $v\in V(G)\setminus \{u\}$.
\end{lemma}
\pf Let $O$ be an orientation of $G$ such that for any orientation $O'$ of $G$,
\[
|{x\in V(G)\setminus \{u\}:d^+_O(x)\notin H(x)}|\leq |{x\in V(G)\setminus \{u\}:d^+_{O'}(x)\notin H(x)}|.
 \]

We claim  ${x\in V(G)\setminus \{u\}:d^+_O(x)\notin H(x)}=\emptyset$. Otherwise, suppose that there exists  $v\in V(G)-\{u\},$ such that  $d^+_O(v)\notin H(v)$.
Then let $P$ be a path connecting $u$ and $v$, and let $P(x,y)$ denote the sub-path of $P$ between two of its vertices $x$ and $y$. Note that we may choose $v$ such that that for any $w\in V(P)\setminus \{u,v\}$, $d^+_O(w)\in H(w)$.
 By (\ref{main-eq}), $d^+_{O\bigtriangleup P}(v)\in H(v)$.	
Let $x$ be the vertex with $d^+_{O\bigtriangleup P}(x)\notin H(x)$ nearest to $v$ on $P$ (if there exists no such vertex, then define $x:=u$). By the choice of $x$, one can see that $d^+_{O\bigtriangleup P(v,x)}(x)\in H(x)$, and for any $y\in V(P(x,v))\setminus \{x,v\}$,  $d^+_{O\bigtriangleup P(x,v)}(y)=d^+_{O\bigtriangleup P}(y)\in H(y)$. Thus
 we have
  \[
  \{w\in V(G)\setminus \{u\}:d^+_{O\bigtriangleup P(x,v)}(w)\notin H(w)\}\subseteq \{y\in V(G)\setminus \{u\}:d^+_O(y)\notin H(y)\}.
   \]
   Note that $v$ is \emph{feasible}. Hence
   \[
   |{y\in V(G)\setminus \{u\}:d^+_{O\bigtriangleup P(x,v)}(y)\notin H(y)}|<| {y\in V(G)\setminus \{u\}:d^+_O(y)\notin H(y)}|,
   \] contradicting to the choice of $O$.  \qed

For $u\in V(G)$, define
\begin{align*}
	D^+_{G,H}(u)=\{d^+_O(u):O\ \mbox{is an orientation of}\ G, \ d^+_O(x)\in H(x)\ if\ x\in V(G)\setminus \{ u\}\},
\end{align*}
i.e.,, $D^+_{G,H}(u)$ is the set of all possible out-degrees in $u$, under the condition that all the other vertices are \emph{feasible}.
 We often use the following  observations:
\begin{obse}\label{Ob1}
\begin{enumerate}
\item  [$($i$)$]  By Lemma \ref{prop1}, $D^+_{G,H}(u)\neq \emptyset.$

\item  [$($ii$)$]  There exists an H-orientation if and only if $D^+_{G,H}(u)\cap H(u)\neq \emptyset$.

\item  [$($iii$)$]  If $i,i+1\in D^+_{G,H}(u)$, then by (\ref{main-eq}), $D^+_{G,H}(u)\cap H(u)\neq \emptyset$, and consequently there exists an \emph{$H$-orientation}.
\end{enumerate}
\end{obse}

Given a \emph{dense pair} $(G,H)$ and $x\in V(G)$, we call $H(x)$ \emph{odd}, or \emph{even}, if $H(x)$ consists of only odd or of only even integers. A vertex $x$ will be said to have \emph{fixed parity} if $H(x)$ is odd or even. When there is no confusion, we also call $x$ is \emph{odd} (or \emph{even}) if $H(x)$ is odd (even, respectively).
\begin{lemma}\label{prop_parity2}
	 Suppose $u\in V(G)$ is not a cut vertex in $G$. Then $D^+_{G,H}(u)$ either contains two consecutive integers, or if not, then it consists of all the odd, or all the even integers in $[0,d_G(u)]$.
\end{lemma}
\pf If $d_G(u)=1$, then the result holds obviously. Next we may assume that $d_G(u)\geq 2$.   Consider $D^+_{G,H}(u)$ contains two consecutive integers. Since $(G,H)$ is dense, $G$ contains an $H$-orientation.

Now suppose that $G$ does not admit $H$-orientation.

\medskip
Claim 1.~$i\in D^+_{G,H}(u),\ i+1\notin D^+_{G,H}(u),\  i+2\leqslant d_G(u)\  \mbox{implies}\  i+2\in D^+_{G,H}(u).$
\medskip

Let $O$ be an orientation of $G$  such that $d_{O}^+(x)\in H(x)$ for all $x\in V(G)\setminus \{u\}$,  $d^+_O(u)=i$, and $e_1,e_2\in \delta(u)$ satisfying that the orientations of them are pointing to $u$ with respect to $O$. Since $u$ is not a cut vertex, there exists a cycle $C$ containing both $e_1$ and $e_2$. If for every $x\in V(C)\setminus \{u\}$, $d_{O\bigtriangleup C}^+(x)\in H(x)$, then by the definition, one can see that $i+2\in D^+_{G,H}(u)$. For $x\notin V(C)\setminus \{u\}$, we use $P(u,x)$ denote the sub-path clockwise along $C$. Let $v\in V(C)\setminus \{u\}$ be the first vertex along $C$ clockwise such that $d_{O\bigtriangleup C}^+(v)\notin H(v)$. Thus   for every  $y\in V(P(u,v))\setminus\{u,v\}$, we have $d_{O\bigtriangleup C}^+(y)=d_{O}^+(y)\in H(y)$. Moreover, since $(G,H)$ is dense, $d_{O\bigtriangleup C}^+(v)\notin H(v)$ and $d_{O}^+(v)\in H(v)$,  we have $d_{O\bigtriangleup P(u,v)}^+(v)\in H(v)$.
Recall that $i\notin H(u)$ and $i+2\leq d_G(u)$.  Since $(G,H)$ is dense, then we have $i+1\in H(u)$. Thus we infer that $d_{O\bigtriangleup P(u,v)}^+(u)\in H(u)$. Now one can see that for all $x\in V(G)$, $d_{O\bigtriangleup P(u,v)}^+(x)\in H(x)$, i.e., $G$ has an $H$-orientation, a contradiction. This completes the proof.

With similar discussion, we have the following statement.

\medskip
Claim 2.~ $i\in D^+_{G,H}(u),\ i-1\notin D^+_{G,H}(u),\  i-2\geqslant 0\  \mbox{implies}\  i-2\in D^+_{G,H}(u).$
\medskip

Combining Claims 1 and 2, the result is followed. This completes the proof.
\qed
\begin{lemma}\label{2-conn-parity}
	If $(G,H)$ is \emph{dense}, and $G$ does not have $H$-orientation, then all the vertices of $G$ except possibly the cut vertices have fixed parity.
\end{lemma}
\pf 
Suppose that $v\in V(G)$ is not a cut vertex.
If $D^+_{G,H}(v)$ contains two consecutive integer, then $D^+_{G,H}(v)\cap H(v)\neq \emptyset$ since $(G,H)$ is dense. So by Observation \ref{Ob1}, we may infer that $G$ has an $H$-orientation, a contradiction.
 If $D^+_{G,H}(v)$ does not contain two consecutive integers, then by Lemma \ref{prop_parity2}, it consists of all the odd or all the even integers of the interval $[0,d_G(v)]$. Since $(G,H)$ is dense and  $H(v)\cap D^+_{G,H}(v)=\emptyset$, we have $H(u)=[0,d_G(v)]\backslash D^+_{G,H}(v)$.
 Thus $H(v)$ has fixed parity. This completes the proof. \qed

\begin{theorem}\label{3}
	Let $G$ is 2-connected and $(G,H)$ is dense. $G$ has no $H$-orientation if and only if $x$ has fixed parity for all $x\in V(G)$, and $|\{x\in V(G)\ |\ H(x) \mbox{ is odd}\}|\neq e(G)\pmod 2$.
\end{theorem}
	
\pf Necessity. Suppose that $G$ has no $H$-orientation. By Lemma \ref{2-conn-parity}, we have that every vertex has fixed parity. Choose an arbitrary $v\in V(G)$. By Lemma \ref{prop1}, there exists an orientation $O$ of $G$ such that $d^+_O(x)$ has the same parity as $H(x)$ for all $x\in V(G)\backslash \{v\}$. If $d^+_O(x)$ and  $H(x)$ has the same parity, then we have $d_O^+(x)\in H(x)$ since $(G,H)$ is dense. Thus we may infer that $G$ has an orientation, a contradiction.

Now consider $d^+_O(x)$ and  $H(x)$ has the different parity.  Then we have
\[
|\{x\in V(G)\ |\ H(x) \mbox{ is odd}\}|\neq \sum\limits_{x\in V(G)}d^+_O(x) \pmod 2.
\]
 Since $\sum\limits_{v\in V(G)}d^+_O(v)=e(G)$, one can see that $|\{x\in V(G)\ |\ H(x) \mbox{ is odd}\}|\neq e(G)\pmod 2$. This completes the proof.

Sufficiency. By contradiction, suppose that $G$ has an $H$-orientation. Then we have $\sum\limits_{v\in V(G)}d^+_O(v)\equiv |\{x\in V(G)\ |\ H(x) \mbox{ is odd}\}|\pmod 2$. Other hand, $\sum\limits_{v\in V(G)}d^+_O(v)=e(G)$. Thus we may infer that $|\{x\in V(G)\ |\ H(x) \mbox{ is odd}\}|\equiv e(G)\pmod 2$, a contradiction. This completes the proof.\qed


We are still far from our goal: in order to see clearly when $G$ does not have \emph{H-orientation}, we should better understand the set  $D^+_{G,H}(v)$ for cut vertices as well.

\section{Connected Graphs}

\noindent \textbf{Proof of Theorem \ref{main}.}
Firstly, we prove necessity. Suppose there exists an $H$-orientation $O$.  Let $(x_1,...,x_k)$ be a parity trace and let $V^1$ be defined as above. Now it is sufficient for us to prove $|V^1|\equiv e(G)\pmod 2$. Since $\sum_{x\in V(G)}d_O^+(x)=e(G)$,  it is enough to prove that $d^+_O(x)$ is odd if and only if $x\in V^1$.

\medskip
\textbf{Claim 1.} Let $V_j^1$ and $V_j^2$ be defined as above. Then $d^+_O(x)$ is odd for $x\in V^1_j$, and $d^+_O(x)$ is even for  $ x\in V^2_j$.
\medskip

 We prove this statement by induction on $j\leqslant k$. For $j=0$, the result holds  by definition. Suppose this is true for $j=i-1$. Let $C$ be a connected component of $G-x_i$ such that $C\subseteq V^1_{i-1} \cup V^2_{i-1}$.

  If  $|C\cap  V^1_{i-1}|+ e(C)\not \equiv d_G(x_i,C)\pmod 2$, then $d^+_O(x_i,C)\geq 1$.  By the definition of ``$l_i$", this immediately implies  $d^+_O(x_i)\geq l_i$.

  If $C\subseteq V^1_{i-1} \cup V^2_{i-1}$, and the parity of $|C\cap  V^1_{i-1}|$ is different from the parity of $e(C)$, then $d^+_O(x_i,C)\leq e_G(x_i,C)-1$. By the definition of ``$t_i$", $ d^+_O(x_i)\leq d_G(x_i)-t_i=u_i$ immediately follows.

Consequently $d^+_O(x_i)\in [l_i,u_i]\cap H(x_i)$, whence $d^+_O(x_i)$ is odd if $x_i\in V^1_i$, and even if $x_i\in V^2_i$. This completes the proof claim 1.
 Since $(x_1,...x_k)$ is a parity trace, $V(G)=V^1\cup V^2$. By Claim 1,   $d^+_O(x)$ is odd if and only if $x\in V^1$. This completes the proof of necessity.

Next we prove sufficiency. By contradiction, suppose that $G$ does not admit  \emph{H-orientation}. It is sufficient to construct a parity trace $(x_1,...x_k)$ satisfying that the parity of $|V^1| \neq e(G)\pmod 2$.
Suppose  for $1\leq i\leq k$,  $V^1_{i-1}$ and $ V^2_{i-1}$ have been defined.
We construct $(x_1,...x_k)$ so that (\ref{parity}) holds. Now we give the definition of $V^1_{i}, V^2_{i}$.

\medskip
\textbf{Claim 2.}	There exists $x_i\in V(G)\setminus (V^1_{i-1}\cup V^2_{i-1})$, for which all components of $G-x_i$ except possibly one are subsets of $V^1_{i-1}\cup V^2_{i-1}$.
\medskip

By Lemma \ref{2-conn-parity}, every vertex in $V(G)\setminus(V_0^1\cup V_0^2)$ is a cut vertex of $G$. Recall that $V_0^1\cup V_0^2\subseteq V_i^1\cup V_i^2$ for $i\geq 0$. Thus all vertices in $V(G)\setminus(V_i^1\cup V_i^2)$ are cut vertices of $G$. Write $U_i=V(G)\setminus(V_i^1\cup V_i^2)$. We choose $x\in U_i$ such that
\[
\max \{|V(C)\cap U_i|\ |\ C \mbox{ is a connected components of $G-x$}\}
\]
 is as large as possible. Let $R$ be the connected component of $G-x$ such that
\[
|V(R)\cap U_i|=\max \{|V(C)\cap U_i|\ |\ C \mbox{ is a connected components of $G-x$}\}.
\]
  We claim $U_i\backslash\{x\}\subseteq V(R)$. Otherwise, let $x'\in U_i\setminus (V(R)\cup\{x\})$. Then $G-x'$ has a connected component $R'$ such that $x\in V(R')$ and $V(R)\subseteq V(R')$, contradicting to the choice of $x$. This completes the proof of claim 2.


Choose now $x_i$ to be a point with the property stated in the Claim 2. Let $\mathcal{B}$ be the set of induced subgraphs with the form $G[V(C)\cup \{x_i\}]$, where $C$ is a connected component of $G- {x_i}$. Thus $x_i$ is not a cut vertex of graph $B$  for every $B\in \mathcal{B}$, and
\begin{align}\label{family}
D^+_{G,H}(x_i)=\sum_{B\in \mathcal{B}}D^+_{B,H}(x_i).
\end{align}
	(If $X,Y,...Z$ are sets of numbers, then $X+Y+...+Z:=\{x+y+...+z\ |\ \ x\in X,\ y\in Y,...,z\in Z\}.)$ Since $G$ does not have \emph{H-orientation},  $D^+_{G,H}(x_i)$ does not contain two consecutive integers. Thus  for all $B\in \mathcal{B}$,  $D^+_{B,H}(x_i)$ also does not contain two consecutive integers.
 By Lemma \ref{prop_parity2}, for each $B\in \mathcal{B}$,  every element of $D^+_{B,H}(x_i)$ is odd or every element of $D^+_{B,H}(x_i)$ is even.  Comparing this with (\ref{family}), and with the definition of $l_i,u_i$, we see that $D^+_{G,H}(x_i)\subseteq [l_i,u_i]$.
	
	On the other hand, by the choice of $x_i$, for all $B\in \mathcal{B}$, $V(B)\setminus \{x_i\}\subseteq V^1_{i-1}\cup V^2_{i-1}$ holds, except for possibly one $B_0\in \mathcal{B}$. Recall that $x$ is not a cut vertex of $B_0$. Then $D^+_{B_0,H}(x_i)$ consists of all the odd or all the even integers in $[0, d_{B_0}(x_i)]$. Thus according to (\ref{family}), $D^+_{G,H}(x_i)$ consists of each second number starting with either $l_i$ or $l_i+1$ and until either $u_i$, or $u_i-1$, depending on the parity of $D^+_{B_0,H}(x_i)$. In other words $D^+_{G,H}(x_i)$ consists of all the odd or all the even integers in $[l_i,u_i]$. Since $G$ does not admit \emph{H-orientation}, we have $D^+_{G,H}(x_i)\cap H(x_i)=\emptyset$. Using (\ref{main-eq}), we have
	\begin{align}\label{reason}
	 [l_i,u_i]\cap H(x_i)=[l_i,u_i]\setminus D^+_{G,H}(x_i),
	\end{align}
	which implies (\ref{parity}) holds.
	
We have proved that if there is no \emph{H-orientation}, then there exists a parity trace $(x_1,...x_k)$. Now it is sufficient to show that for this parity trace, $|V^1|\neq e(G)\pmod 2$.
	
	Without loss generality, we may assume that for every $x\in V(G)\setminus x_k =V^1_{k-1} \cup V^2_{k-1}$, $d_O^+(x)\in H(x)$. Since $(x_1,...x_k)$ is a parity trace, applying Claim 1, we know that
 for $x\in V^1_{k-1}$, $d^+_O(x)$ is odd and for $x\in V^2_{k-1}$, $d^+_O(x)$ is even. Since $G$ does not admit $H$-orientation, then $d^+_O(x_k)$ has the different parity with $H(x_k)\cap [l_k,u_k]$. Thus we may infer that $\sum\limits_{v\in V(G)}d^+_O(v)\neq |V^1|\pmod 2$. Recall that
 $\sum\limits_{v\in V(G)}d^+_O(v)=e(G)$. Hence we get $|V^1|\neq e(G)\pmod 2$, this completes the proof. \qed


 We can derive the following result, which  partly confirms Conjecture \ref{AK}.
 \begin{coro}\label{cor1}
 		Let $G$ be a graph and let $F,H:V(G)\to 2^N$ such that $(G,H)$ is \emph{dense} and   $F(v)=[0,d_G(v)]\setminus H(v)$ for all $v\in V(G)$. If for every $v\in V(G)$,
 	$$|F(v)|\leq \frac{1}{2}(d_G(v)-1), $$
 	then $G$ admits an $H$-orientation.
 \end{coro}
\pf  Since $|F(v)|\leq \frac{1}{2}(d_G(v)-1)$, then $H(v)$ contains two consecutive integers.
 Let $v_0\in V(G)$ and $v_0$ is not a cut vertex of $G$. Since $H(v_0)$ contains two consecutive integers, then  by Lemma \ref{2-conn-parity}, $G$ admits an $H$-orientation. \qed 

\medskip

\noindent \textbf{Remark.} From the proof of Corollary \ref{cor1}, the result still holds if the following conditions holds.
\begin{itemize}
  \item  For all $x\in V(G)$,  $|F(x)|\leq \frac{1}{2}(d_G(v)+1)$ and $F(x)$ contains no two consecutive integers;

  \item there exists
$u\in V(G)$ such that $u$ is not cut vertex, $|F(u|\leq \frac{1}{2}(d_G(u)-1)$.
\end{itemize}

By Theorem \ref{main}, we may obtain the following result.
\begin{coro}\label{cor2}
 		Let $G$ be a graph and let $H:V(G)\to 2^N$ such that for every $v\in V(G)$,  $H(v)\in  \{\{0,2,...,2\lceil d_G(v)/2\rceil\}, \{1,3,\ldots, 2\lceil d_G(v)/2\rceil-1\}\} $. If $|\{v\in V(G)\ |\ H(v) \mbox{ is odd}\}|\equiv e(G)\pmod 2$, then $G$ admits an $H$-orientation.
 \end{coro}

\end{document}